\documentclass[12p]{article}
\usepackage[latin1]{inputenc}
\usepackage{a4wide}
\usepackage{amsmath}
\usepackage{amssymb}
\usepackage{amsfonts}
\usepackage{amsthm}
\usepackage{url}
\usepackage{graphicx}

\newcommand{\E}{\ensuremath{\mathbb{E}}}
\newcommand{\V}{\ensuremath{\mathrm{Var}}}
\newcommand{\Prob}{\ensuremath{\mathbb{P}}}

\begin{document}
\newtheorem{theorem}{Theorem}[section]
\newtheorem{korollar}[theorem]{Corollary}
\newtheorem{proposition}[theorem]{Proposition}

\title{\bf Asymptotic analysis of Hoppe trees}
\author{Kevin Leckey and Ralph Neininger\thanks{Email: \{leckey, neiningr\}@math.uni-frankfurt.de}\\
Institute for Mathematics\\
J.W.~Goethe University Frankfurt\\
60054 Frankfurt am Main\\
Germany} 

\date{July 5, 2012}
\maketitle

\begin{abstract}
We introduce and analyze a random tree model associated to Hoppe's urn. The tree is built successively by adding nodes to the existing tree when starting with the single root node. In each step a node is added to the tree as a child of an existing node where these parent nodes are chosen randomly with probabilities proportional to their weights. The root node has weight $\vartheta>0$, a given fixed parameter, all other nodes have weight $1$. This resembles the stochastic dynamic of Hoppe's urn. For $\vartheta=1$ the resulting tree is the well-studied random recursive tree. 
We analyze the height, internal path length and number of leaves of the Hoppe tree with $n$ nodes as well as the depth of the last inserted node asymptotically as $n\to \infty$. Mainly  expectations, variances and asymptotic distributions of these parameters are derived.
\end{abstract}

\noindent
{\em  AMS 2010 subject classifications.} Primary 60F05, 60C05; secondary  60G42, 68R05.\\
{\em Key words.} Hoppe urn, random tree, weak convergence, martingale, combinatorial probability.

\section{Introduction}
We consider a random tree model associated and derived from Hoppe's urn: In
Hoppe's urn, see \cite{Ho86},  there initially is one red ball. In each step one of the balls is   drawn from the urn independently with probabilities proportional to the weights of the balls. The red ball has weight $\vartheta>0$, all other balls have weight $1$. Here the parameter $\vartheta>0$ is given and fixed throughout the evolution of the urn. When a ball is drawn it is placed back to the urn together with a ball of the same color unless the ball drawn is the red ball. In this case the red ball is placed back together with a ball of a new color not yet being present in the urn.
This model has been introduced for deriving and interpreting the Ewens sampling formula and is related to the infinite alleles model in population genetics, the parameter  $\vartheta>0$ modeling the mutation rate. The decomposition of the balls into  groups of the same color (neglecting the red ball) leads to a Chinese restaurant process, the $(0,\vartheta)$ seating plan, see Pitman \cite[page 61]{Pi06}.

A random tree model, which we subsequently call Hoppe tree, is associated to the Hoppe urn as follows: The balls in the urn are represented by nodes in the tree. Each node $v$  is child of node $w$ in the tree if the ball corresponding to $v$ was placed first in the urn together with the ball corresponding to $w$ when the $w$-ball was drawn. In other words the tree grows successively: In each step a node is chosen independently and with probability proportional to the weights of the nodes (the root having weight $\vartheta$, all other nodes having weight $1$) and a new node is added as child of the chosen node. For $\vartheta=1$ this is a well-known and well-studied random tree model, the random recursive tree, see, e.g., Smythe and Mahmoud \cite{MS94}.

The aim of the present note, which is based on the first author's master's thesis \cite{Le11}, is to study asymptotic properties of the Hoppe tree as its size $n$ tends to infinity. In particular we are interested in the deviation from the random recursive tree model caused by the perturbation of the root weight from $\vartheta=1$ to $\vartheta\neq 1$. As characteristics of the tree we study the depth $D_n^{(\vartheta)}$ of the $n$-th inserted node in the  tree, defined as its distance to the root of the tree. Furthermore the tree's height $H_n^{(\vartheta)}$ is studied, which is the maximal depth $\max_{1\le i\le n} D_i^{(\vartheta)}$, its internal path length $I_n^{(\vartheta)}=\sum_{1\le i\le n} D_i^{(\vartheta)}$ and the number of leaves of the tree. A node is a leaf if it has no child in the tree. Our results show, that the perturbation of the root weight does typically not affect the first order behavior of the quantities, an exception being the variance and limit law of the internal path length. Hence, we give second order expansions to reveal the asymptotic dependence on $\vartheta$. 

The paper is organized as follows: In the second section the results on the four quantities mentioned above are stated, the proofs being collected in the third section. 

\subsubsection*{Acknowledgment} We thank Henning Sulzbach for comments on a draft of this note and two anonymous referees for their careful reading.

\section{Results}
In this section the results on depth, height, internal path length and number of leaves are stated. 
Throughout the parameter $\vartheta>0$ is arbitrary and fixed.
All asymptotic statements as well as the use of the Bachmann-Landau symbols are understood as $n$, the number of nodes in the Hoppe tree, tends to infinity. Moreover, we use the digamma and trigamma functions $\Psi=\frac{d}{dx}\log \Gamma$ and $\Psi_1=\frac{d^2}{dx^2}\log \Gamma$ respectively. By the properties of the digamma and trigamma functions, see e.g. \cite[6.3. and 6.4.]{AS64}, we have
\begin{align*}
\sum_{i=1}^{n-2} \frac 1 {\vartheta+i} &= \Psi (\vartheta+n-1)-\Psi(\vartheta+1) = \log n - \Psi (\vartheta +1)+o(1),\\
\sum_{k=1}^\infty \left(\frac 1 {\vartheta+k}\right)^{2} &=  \Psi ' (\vartheta+1) = \Psi_1 (\vartheta+1).
\end{align*}
\subsubsection*{Depth of a node}
For the depth $D_n^{(\vartheta)}$ we have a  distributional representation as sum of independent Bernoulli variables: 
\begin{theorem}
\label{d1}
For the depth $D_n^{(\vartheta)}$  of the $n$-th node in a Hoppe tree we have for all $n\geq 2$ 
$$D_n^{(\vartheta)}\stackrel{d}{=} 1 + \sum_{i=1}^{n-2} B_{i},$$
where $B_1,\ldots B_{n-2}$ are independent and $\Prob(B_i=1)=1-\Prob(B_i=0)= \frac 1 {\vartheta+i}$ for $i=1,\ldots,n$.
\end{theorem}
Asymptotic results can hence easily be obtained, e.g., the following. We denote by $\Pi(\lambda)$  the Poisson distribution with parameter $\lambda>0$, by $d_\mathrm{TV}$ the total variation distance between probability measures, by $\stackrel{d}{\longrightarrow}$ convergence in distribution and by $\mathcal{N}(0,1)$ a real random variable with the standard normal distribution. 
\begin{korollar}
\label{d2}
The depth $D_n^{(\vartheta)}$ of the $n$-th node in a Hoppe tree satisfies
\begin{align}
\E[D_n^{(\vartheta)}] &= 1 + \sum_{i=1}^{n-2} \frac 1 {\vartheta+i} = \log n -\Psi(\vartheta+1) +1 + o(1), \notag \\
\V(D_n^{(\vartheta)}) &= \sum_{i=1}^{n-2}\frac 1 {\vartheta+i} - \sum_{i=1}^{n-2} \left( \frac 1 {\vartheta+i} \right)^2 \notag \\
&= \log n -\Psi(\vartheta+1)-\Psi_1(\vartheta+1) + o(1), \notag \\
\frac {D_n^{(\vartheta)}-\E[D_n^{(\vartheta)}]} { \sqrt{\V(D_n^{(\vartheta)})}} & \stackrel{d}{\longrightarrow}  \mathcal{N}(0,1), \label{CLTD} \\
d_\mathrm{TV} \left(\mathcal{L} (D_n^{(\vartheta)}), \Pi \left(\E[D_n^{(\vartheta)}]\right) \right) &= \mathcal{O} \left(\frac 1 {\log n } \right). \notag
\end{align}
\end{korollar}

\subsubsection*{Height of the Hoppe tree}
The height $H_n^{(\vartheta)}$ of the Hoppe tree can be analyzed by drawing back to results on the height for random recursive trees, see Addario-Berry and Ford \cite{AF10}, in particular they show that 
\begin{align}
M_n:=\E[H_n^{(1)}]= e \log n - \frac 3 2 \log\log n + \mathcal{O}(1) \label{ERRT}
\end{align}
as $n\to\infty$. We transfer their results to arbitrary $\vartheta>0$:
\begin{theorem}
\label{h1}
For the height $H_n^{(\vartheta)}$  of a Hoppe tree with $n$ nodes we have: For all $\alpha<\frac {1} {3e}$, $\beta < \frac 1 {2e}$  there exist constants $C_\alpha, C_\beta>0$ such that for all $t> 0$
\begin{align*}
\Prob\left(H_n^{(\vartheta)}-M_n\geq t \right) \leq C_\beta e^{-\beta t},\qquad 
\Prob\left(H_n^{(\vartheta)}-M_n\leq -t \right) \leq C_\alpha e^{-\alpha t}.
\end{align*}
The constant $C_\beta$ can be chosen independently of $\vartheta$.
\end{theorem}
\begin{korollar}
\label{h2}
The height $H_n^{(\vartheta)}$ of a Hoppe tree with $n$ nodes satisfies
\begin{align*}
\E[H_n^{(\vartheta)}]= e \log n - \frac 3 2 \log\log n +\mathcal{O}(1), \qquad
\V(H_n^{(\vartheta)}) = \mathcal{O}(1).
\end{align*}
\end{korollar}

\subsubsection*{Number of leaves}
The number of leaves in a Hoppe tree is related to a two-color urn model. We obtain:
\begin{theorem}
\label{l1}
Let $L_n^{(\vartheta)}$ be the number of leaves in a Hoppe tree with $n\geq 2$ nodes. Then
\begin{align}
\E[L_n^{(\vartheta)}] &= \frac {n} 2+ \frac {\vartheta-1} 2 + \mathcal{O} \left( \frac 1 n \right), \nonumber\\
\V( L_n^{(\vartheta)} ) &= \frac {n} {12}+ \frac {\vartheta-1} {12} + \mathcal{O} \left( \frac 1 n \right), \nonumber\\
\Prob(|L_n- \E[L_n]| \geq t )  &\leq  2 \exp \left( - \frac {6 t^2}{n+\vartheta+1} \right) \text{ for all } t>0, n\ge 1, \label{tail_bd} \\
\frac {L_n^{(\vartheta)}-\E[L_n^{(\vartheta)}]} {\sqrt{\V(L_n^{(\vartheta)})}} &\stackrel {d} {\longrightarrow}  \mathcal{N}(0,1). \nonumber
\end{align}
\end{theorem}

\subsubsection*{Internal path length}
Moments of the internal path length can be obtained from our results on the depths of nodes. 
\begin{theorem}
\label{i1}
The internal path length $I_n^{(\vartheta)}$ of a Hoppe tree with $n$ nodes satisfies
\begin{align*}
\E[I_n^{(\vartheta)}]&= (\vartheta+n-1)\sum_{i=1}^{n-1} \frac 1 {\vartheta+i} = n \log n -\Psi(\vartheta+1) n + o(n), \\
\V(I_n^{(\vartheta)}) &= \left( \frac 2 {\vartheta+1} - \Psi_1(\vartheta+1)\right) n^2 + o(n^2).
\end{align*}
Moreover, 
\begin{align*}
\left( \frac {I_n^{(\vartheta)}-\E[I_n^{(\vartheta)}]} {\vartheta+n-1} \right)_{n\geq 1}
\end{align*} is a zero-mean martingale.
\end{theorem}

The internal path length can be analyzed either via martingale methods or the recursive distributional decomposition explained in Figure 1 which allows to apply the contraction method.
\begin{theorem} \label{llipl}
The internal path length $I_n^{(\vartheta)}$ of a Hoppe tree with $n$ nodes satisfies
\begin{align*}
\frac{I_n^{(\vartheta)}-n\log n} {n} \to X^{(\vartheta)}
\end{align*}
for a non-degenerate random variable $X^{(\vartheta)}$, where the convergence holds 
 almost surely and in $L_2$. The distribution ${\cal L}(X^{(\vartheta)})$ is the only integrable solution of the distributional fixed point equation
 \begin{align}
\label{rec_limit_ipl}
 X^{(\vartheta)} \stackrel{d}{=} (1-B)X^{(\vartheta)}+ B \widetilde{X}^{(1)} + B\log(B)+(1-B)\log(1-B)+B,
\end{align}
where $X^{(\vartheta)},\widetilde{X}^{(1)}$ and $B$ are independent, $B$ has the beta$(1,\vartheta)$ distribution and $\widetilde{X}^{(1)}$ is distributed as $X^{(1)}$. For $\vartheta\neq 1$, the solution of (\ref{rec_limit_ipl}) is even unique without integrability assumption.
 \end{theorem}

 \begin{theorem}\label{density}
The limit distribution ${\cal L}(X^{(\vartheta)})$ in Theorem \ref{llipl} has a  Lebesgue density $f_\vartheta$, which is in the Schwartz space on $\mathbb{R}$, i.e., $f_\vartheta$ is infinitely differentiable and together with all its derivatives rapidly decreasing.
\end{theorem}

\section{Proofs}
In the analysis of the tree below the random decomposition of the Hoppe tree shown in Figure 1  is used: The tree is decomposed into the subtree of the second inserted node (left dashed box) and the remaining part of the tree (right dashed box). The stochastic dynamic of the Hoppe tree with parameter $\vartheta$ implies that conditioned on the size $N_n$ of the subtree of the second inserted node this subtree is a random recursive tree, whereas the remaining part is a Hoppe tree with parameter $\vartheta$ and size $n-N_n$. Moreover, conditional on $N_n$ these two trees are independent.
We have the asymptotic behavior
\begin{align}\label{lim_Nn}
\frac{N_n}{n} \to B \; \mbox{almost surely} \quad (n\to \infty)
\end{align} 
where $B$ has the beta$(1,\vartheta)$ distribution having Lebesgue density $x\mapsto \vartheta (1-x)^{\vartheta-1}$, $x\in[0,1]$, see Donnelly and Tavar{\'e} \cite{DT86}.

\begin{figure}\label{fig}
	    \centering
	    \includegraphics[width=9cm]{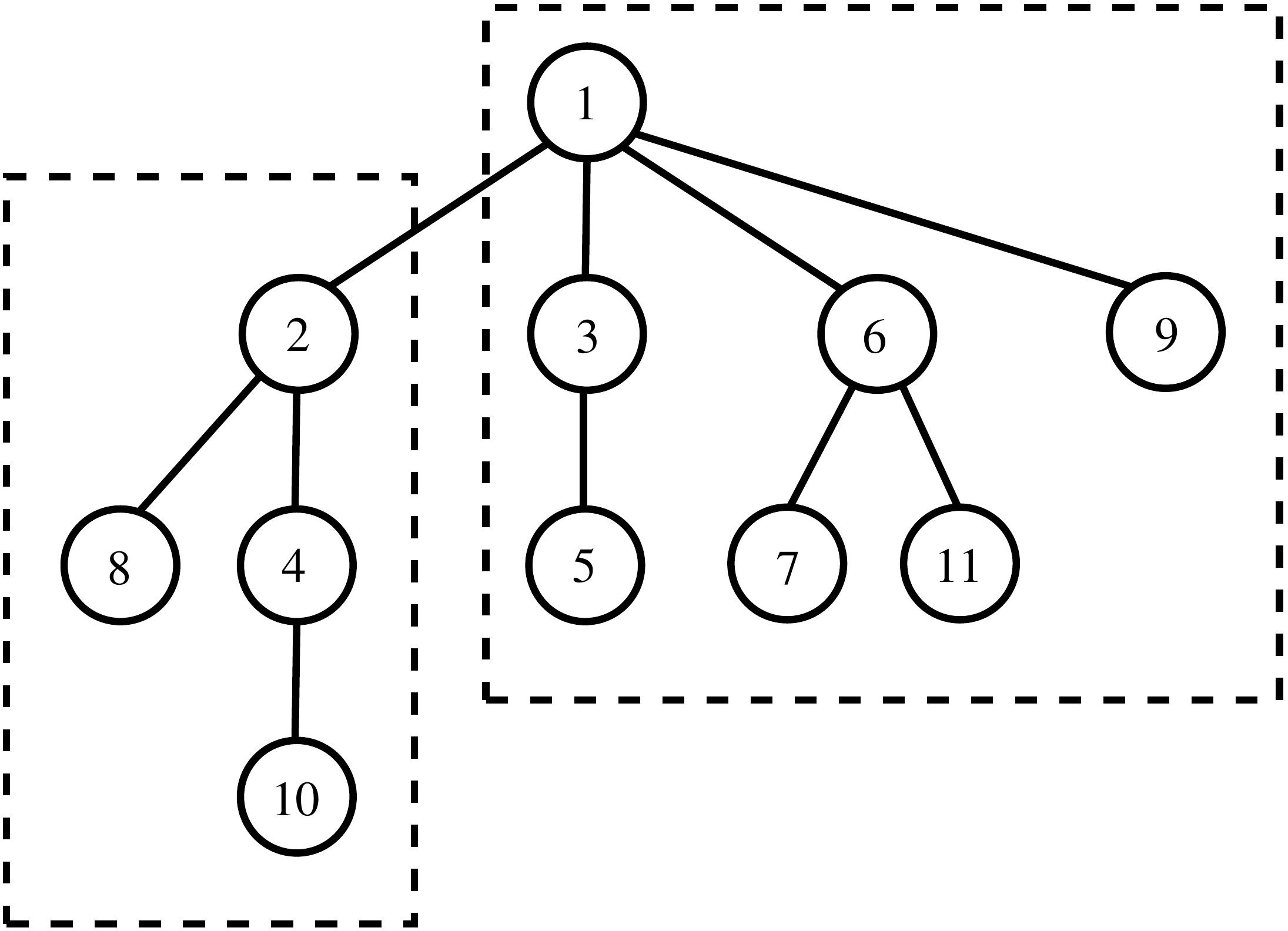}
	    \caption{A Hoppe tree with $11$ nodes. The decomposition into the subtree rooted at node labelled $2$ and the remaining part of the tree is indicated in dashed boxes.}
\end{figure}

\begin{proof}[\textbf{Proof of Theorem \ref{d1}}] We calculate the depth of a node by counting its ancestors in the tree. We have $D_n^{(\vartheta)}= \sum_{i=1}^{n-1} \textbf{1}_{A_{i,n}}$, where $A_{i,j}$ denotes the event that node $i$ is an ancestor of node $j$, $i<j$. Cleary $\Prob(A_{1,n})=1$. Moreover, $\Prob(A_{i,i+1})= \frac 1 {\vartheta + i -1 }$ for $i\geq 2$ by definition of the Hoppe tree. For general $i<n$ let $\xi_{i,n}$ be the number of descendants of node $i$ in a Hoppe tree with $n$ nodes, i.e. the size of subtree rooted in $i$ minus 1. By the dynamics of the Hoppe tree we have
\begin{equation} \label{PAin} \Prob(A_{i,n}|\xi_{i,n-1})=\frac {1+\xi_{i,n-1}} {\vartheta + n-2}. \end{equation}
We calculate $\E[\xi_{i,n-1}]$ by the recursion
$$ \E[\xi_{i,n-1}]= \E[\xi_{i,n-2}+\textbf{1}_{A_{i,n-1}}] = \E[\xi_{i,n-2}]+ \frac {1+\E[\xi_{i,n-2}]} {\vartheta+n-3}.$$
This yields $ \E[\xi_{i,n-1}]=\frac{\vartheta + n - 2}{\vartheta +i-1}-1$ and therefore, by equation (\ref{PAin}),
\begin{equation} \label{PAin2} \Prob(A_{i,n})=\frac 1 {\vartheta+i-1}. \end{equation}
It remains to show that $A_{2,n},\ldots,A_{n-1,n}$ are independent. Note that for $i<j$, $A_{i,j}$ only depends on where the nodes $i+1,\ldots,j$ are inserted. Therefore, we get for all $2\leq k \leq n-2$ and $2\leq i_1 <\ldots <i_k \leq n-1$ independence of $A_{i_1, i_2}, A_{i_2,i_3}, \ldots , A_{i_k , n}$. Since $\bigcap_{j=1}^k A_{i_j, n}$ occurs if and only if $i_j$ is an ancestor of $i_{j+1}$ for every $j\leq k-1$ and $i_k$ is an ancestor of $n$ we have
\begin{align*}
\Prob\left( \bigcap_{j=1}^k A_{i_j, n} \right) &= \Prob\left( A_{i_1, i_2} \cap A_{i_2,i_3}\cap \ldots \cap A_{i_k , n} \right) \\
&= \Prob( A_{i_1, i_2}) \cdot \Prob(A_{i_2,i_3}) \cdot \ldots \cdot \Prob( A_{i_k , n}) \\
&= \prod_{j=1}^k \Prob(A_{i_j,n}),
\end{align*}
where (\ref{PAin2}) is used in the last equation.
With $B_i= \textbf{1}_{A_{i+1,n}}$ and $\textbf{1}_{A_{1,n}}=1$ this yields the assertion. 

For related reasoning in the analysis of the depth  in other random tree models see Dobrow and Smythe \cite{DoSm95}.
\end{proof}

\begin{proof}[\textbf{Proof of Corollary \ref{d2}}] Theorem \ref{d1} implies expectation and variance of $D_n^{(\vartheta)}$. Moreover, by Lindeberg's version of the central limit theorem  (CLT) we obtain the CLT for $D_n^{(\vartheta)}$ in (\ref{CLTD}) and by \cite[Equation (1.23)]{BHJ92} we get $d_\mathrm{TV} (\mathcal{L} (D_n^{(\vartheta)}), \Pi ( \E[D_n^{(\vartheta)}]) ) = \mathcal{O}( 1 / \log n )$.
\end{proof}

\begin{proof}[\textbf{Proof of Theorem \ref{h1}}] Addario-Berry and Ford showed in \cite[Corollary 1.3]{AF10} that the expected height $M_n:=\E[H^{(1)}]$ of a random recursive tree satisfies (\ref{ERRT}) and that for all $c'< \frac 1 {2e}$ there exists a constant $C=C(c')$ such that for all $n\geq 1 $ and $t>0$
$$\Prob(|H_n^{(1)}- M_n|\geq t)\leq C e^{-c' t}.$$
Recall that in a Hoppe tree with $n\geq 1$ nodes and parameter $\vartheta>0$ by $N_n$  the size of the subtree rooted in node 2 is denoted and that this subtree, conditioned on its size, is a random recursive tree.  

By an obvious coupling argument between Hoppe trees for different parameters $\vartheta$ we have $H_n^{(\vartheta_1)}\preccurlyeq H_n^{(\vartheta_2)}$ for all $\vartheta_1 \geq \vartheta_2$, where $\preccurlyeq$ denotes stochastic domination. In the extremal case $\vartheta=0$ (for definition of the tree start with the root and one child) we obtain $H_n^{(\vartheta)}\preccurlyeq H_n^{(0)} \stackrel{d}{=}1+H_{n-1}^{(1)} \preccurlyeq 1+ H_n^{(1)}$. Therefore, we get $\Prob\left(H_n^{(\vartheta)}-M_n\geq t \right) \leq \widehat{C} e^{-c' t}$, $\widehat{C}= C e^{c'}$,  using the result for random recursive trees. 

In order to prove the left tail inequality let $H_{N_n}^{(1)}$ be the height of the subtree rooted in node 2. From $H_n^{(\vartheta)}\geq H_{N_n}^{(1)}$ we obtain for all $t>0$ and $\alpha>0$ (later we have to restrict to $\alpha$ as in the Theorem)
\begin{eqnarray*}
\Prob(H_n^{(\vartheta)}-M_n \leq -t)&\leq&\Prob(\{H_{N_n}^{(1)}-M_n\leq-t \} \cap \{ N_n\geq e^{-\alpha t} n \})\\
 & & +\Prob(\{H_{N_n}^{(1)}-M_n \leq -t \} \cap \{ N_n< e^{-\alpha t} n \}),\\
 &\leq & \Prob(H_{\lceil e^{-\alpha t} n \rceil}^{(1)}-M_n \leq -t)+\Prob(N_n< e^{-\alpha t} n).
\end{eqnarray*}
 Again, by using the result for random recursive trees and $M_n- \E[H_{\lceil e^{-\alpha t} n \rceil}^{(1)}]=e \alpha t + \mathcal{O}(1)$ we obtain for  $\alpha=c'/(1+ec')$  a constant $C_1$ such that
$$ \Prob(H_{\lceil e^{-\alpha t} n \rceil}^{(1)}-M_n \leq -t) \leq C_1 e^{-c'(1-e\alpha) t}= C_1 e^{-\alpha t}.$$
Hence we have such an upper bound for all $\alpha<1/(3e)$. 
To get an upper bound for $\Prob(N_n< e^{-\alpha t} n)$ note that for all $1 \leq k\leq n-1$ 
$$\Prob(N_n=k) = \binom {n-2} {k-1} \frac { \vartheta (\vartheta+1) \cdots (\vartheta + n -(k+2) ) (k-1)!} { (\vartheta +1) \cdots (\vartheta +n-2)}.$$
This yields for all $\varepsilon\in (0,1)$ that
$$\Prob(N_n \leq \varepsilon n ) \leq 3(\vartheta+1) \varepsilon .$$
Therefore,
$$\Prob(H_n^{(\vartheta)}-M_n \leq -t) \leq (C_1+3(\vartheta+1)) e^{-\alpha t}.$$
This implies the assertion.
\end{proof}

\begin{proof}[\textbf{Proof of Corollary \ref{h2}}] By Theorem \ref{h1} we have
$$ \E[|H_n^{(\vartheta)}-M_n |]=\mathcal{O}(1).$$
Consequently, $ \E[H_n^{(\vartheta)}]=M_n+\mathcal{O}(1)=e \log n - \frac 3 2 \log \log n + \mathcal{O}(1)$. \\
Moreover, the tail bound from Theorem \ref{h1} implies
$$\V (H_n^{(\vartheta)})\leq \E[(H_n^{(\vartheta)}-M_n)^2] = \mathcal{O}(1).$$

\end{proof}
For the proof of the tail bound in Theorem \ref{l1} we use the following version of Azuma-Hoeffding's inequality with conditional ranges:
\begin{proposition}
 \label{azumahoef}
Let $W_1,\ldots,W_n$ be a  martingal difference sequence with respect to a filtration $(\mathcal{F}_i)_{0\leq i \leq n}$ with $\mathcal{F}_0 =\{\emptyset , \Omega \}$. Suppose that for every $1\leq i \leq n$ there exists a constant $c_i\geq 0$ and an $\mathcal{F}_{i-1}$ measurable random variable $Z_i$ such that $Z_i \leq W_i \leq Z_i + c_i$ almost surely. Then we have for all $t>0$
$$\Prob\left( \left| \sum_{i=1}^n W_i \right| \geq t \right) \leq 2 \exp \left(-\frac {2t^2}{\sum_{i=1}^n c_i^2}\right).$$
\end{proposition}

\begin{proof}[\textbf{Proof of Theorem \ref{l1}}] We have $L_n^{(\vartheta)}=L_{n-1}^{(\vartheta)}+Y_n$, where
$$Y_n=\begin{cases} 1, & \text{if the parent of node $n$ was not a leaf at time $n-1$,} \\ 0, & \text{otherwise.} \end{cases}$$
Therefore, for $n\geq 2$, almost surely
$$ \E[L_{n+1}^{(\vartheta)}|L_1^{(\vartheta)},\ldots, L_{n}^{(\vartheta)}]= L_{n}^{(\vartheta)}+ 1-\frac {L_n^{(\vartheta)}}{\vartheta + n -1}=\frac{\vartheta+n-2}{\vartheta+n-1} L_n^{(\vartheta)}+1.$$
With
\begin{align} \label{scal_mg}
X_n=(\vartheta+n-2)\left(L_n^{(\vartheta)}-\left(\frac {n-1} 2 +\frac {\vartheta(n-1)}{2(\vartheta+n-2)}\right)\right)
\end{align}
the sequence $\left( X_n \right)_{n\geq 2}$ is a zero-mean martingale and 
$$ \E[L_n^{(\vartheta)}]=\frac {n-1} 2 +\frac {\vartheta(n-1)}{2(\vartheta+n-2)}= \frac {\vartheta+n-1} 2 + \mathcal{O}\left(\frac 1 n \right).$$
With the representation 
$$X_i-X_{i-1}=(\vartheta+i-2)(Y_i-\E[Y_i])+L_{i-1}^{(\vartheta)}-\E[L_{i-1}^{(\vartheta)}], \quad i\ge 3$$
we have $Z_i \leq X_i-X_{i-1} \leq Z_i + \vartheta+i-2$ where $Z_i=L_{i-1}^{(\vartheta)}-\E[L_{i-1}^{(\vartheta)}]-(\vartheta+i-2)\E[Y_i]$. By Proposition \ref{azumahoef} we have for all $t>0$
$$\Prob(|X_n |\geq t) \leq 2 \exp \left( - \frac {2 t^2}{\sum_{i=3}^n (i+\vartheta-2)^2} \right).$$
Using that the sum in the denominator of the latter exponent is bounded by $(n+\vartheta-2)^3/3+(n+\vartheta-2)^2$ and the scaling in (\ref{scal_mg}) this implies the bound (\ref{tail_bd}).

In order to compute $\V (L_n^{(\vartheta)})$ we have $X_n= \frac {\vartheta+n-2}{\vartheta+n-3} X_{n-1}+ (\vartheta+n-2) (Y_n- \E[Y_n])$. Hence,
\begin{align}
\label{Xn2}
\E[X_n^2]=&\left(\frac {\vartheta+n-2}{\vartheta+n-3}\right)^2 \E[ X_{n-1}^2]+2\frac {(\vartheta+n-2)^2}{\vartheta+n-3} \E[ X_{n-1}(Y_n-\E[Y_n])] \nonumber\\
&~+(\vartheta+n-2)^2 \V (Y_n).
\end{align}
Using $ \E[X_{n-1}]=0$ we have
\begin{align*}
\E[ X_{n-1}(Y_n-\E[Y_n])] &= \E[X_{n-1} \E[Y_n | L_1^{(\vartheta)},\ldots ,L_{n-1}^{(\vartheta)}]]
= \E\left[X_{n-1} \left(1-\frac {L_{n-1}^{(\vartheta)}}{\vartheta+n-2}\right)\right] \\
&=~-\frac 1 {(\vartheta+n-2)(\vartheta+n-3)} \E[X_{n-1}^2].
\end{align*}
Moreover, $ \E[Y_n]=1-\frac {\E[L_{n-1}^{(\vartheta)}]}{\vartheta+n-2}=\frac 1 2 + \mathcal{O}\left( 1 /{n^2}\right)$ and $\V (Y_n)= \frac 1 4 + \mathcal{O} \left(  1 / {n^2} \right)$.\\
Solving (\ref{Xn2}) by the substitution $Q_n= \frac {\vartheta+n-3}{\vartheta+n-2} \E[X_n^2]$ yields
$$\V (L_n^{(\vartheta)})=\frac { \vartheta + n- 1}{12} +\mathcal{O} \left(\frac 1 n \right).$$
To obtain the CLT for $L_n^{(\vartheta)}$ the representation
$$\frac {L_n^{(\vartheta)}-\E[L_n^{(\vartheta)}]} { \sqrt{\V (L_n^{(\vartheta)})}} = \frac {X_n} {\sqrt{\V (X_n)}}$$
allows to apply a general martingale CLT, see, e.g., Hall and Heyde \cite[Theorem 3.2]{HH80}. It is sufficient to show that
$$\Delta_{n,i}:=\frac 1 {\sqrt{\V (X_n)}} (X_i-X_{i-1}), \qquad n\geq 3, 3\leq i \leq n,$$
satisfies
\begin{align*}
\mbox{(a)} \;\; \max_{3\le i \le n} |\Delta_{n,i}| \stackrel {\Prob} {\longrightarrow} 0, \quad\quad
\mbox{(b)}   \;\; \sum_{3\le i \le n} \Delta_{n,i}^2 \stackrel {\Prob} {\longrightarrow} 1,\quad\quad
\mbox{(c)}   \;\;  \max_{n\ge 3}\; \E \!\left[\max_{3\le i \le n} \Delta_{n,i}^2 \right]<\infty.
\end{align*}
For (a) and (c) we have $|X_i-X_{i-1}|=|L_i^{(\vartheta)}-\E[L_i^{(\vartheta)}]+ (\vartheta+i-3)(Y_i-\E[Y_i])|\leq \vartheta +2n+3$ for $i\leq n$ and $\V (X_n) = (\vartheta +n -1)^2 \V (L_n^{(\vartheta)}) \sim \frac {n^3}{12}$. Hence, $|\Delta_{n,i}| \leq  {( 2n +\vartheta+3)}/{\sqrt{\V (X_n)}}$ a.s., which yields that $\max_i |\Delta_{n,i}| \stackrel {\Prob} {\rightarrow} 0$ and that $\E\left[\max_i \Delta_{n,i}^2 \right]$ is bounded in $n$. \\
To compute $\sum_i \Delta_{n,i}^2$ note that by (\ref{tail_bd}) and the Borel-Cantelli Lemma we have $(L_n^{(\vartheta)}-\E[L_n^{(\vartheta)}])/{n} \rightarrow 0 $ almost surely. Hence, for all $n\geq 3$,
\begin{align}
\sum_{i=3}^n \Delta_{n,i}^2 =& \frac 1 {\V (X_n)} \sum_{i=3}^n (L_i^{(\vartheta)}-\E[L_i^{(\vartheta)}])^2 + \frac 2 {\V (X_n)} \sum_{i=3}^n (L_i^{(\vartheta)}-\E[L_i^{(\vartheta)}])(\vartheta+i-3)(Y_i-\E[Y_i]) \notag \\
& + \frac 1 {\V (X_n)} \sum_{i=3}^n (\vartheta+i-3)^2 (Y_i-\E[Y_i])^2. \label{delta2}
\end{align}
By $(L_n^{(\vartheta)}-\E[L_n^{(\vartheta)}])/{n} {\rightarrow} 0 $, $\V (X_n) \sim \frac {n^3}{12}$ and the Cesàro mean we have for the first summand in (\ref{delta2})
$$ \frac 1 {\V (X_n)} \sum_{i=3}^n (L_i^{(\vartheta)}-\E[L_i^{(\vartheta)}])^2 \leq \frac{n^3}{\V (X_n)} \frac 1 n \sum_{i=3}^n \left( \frac {L_i^{(\vartheta)}-\E[L_i^{(\vartheta)}]} i \right)^2 \rightarrow 0,$$
and for the second summand in (\ref{delta2})
\begin{align*}
\lefteqn{\left|\frac 2 {\V (X_n)} \sum_{i=3}^n (L_i^{(\vartheta)}-\E[L_i^{(\vartheta)}])(\vartheta+i-3)(Y_i-\E[Y_i]) \right|}\\
&\leq \frac {2 n^2 (\vartheta+n+3)}{\V (X_n)} \frac 1 n \sum_{i=3}^n \left| \frac {L_i^{(\vartheta)}-\E[L_i^{(\vartheta)}]} i \right| \rightarrow 0.
\end{align*}
Because $\E[Y_i]= \frac 1 2 + \mathcal {O} \left (\frac 1 {i^2} \right)$ we have $(Y_i-\E[Y_i])^2= \frac 1 4 + \mathcal{O} \left( \frac 1 {i^2} \right)$ a.s. and therefore for the last summand in (\ref{delta2}), a.s.
$$\frac 1 {\V (X_n)} \sum_{i=3}^n (\vartheta+i-3)^2 (Y_i-\E[Y_i])^2 {\to} 1 .$$
This implies $\sum_i \Delta_{n,i}^2 \stackrel {\Prob} {\longrightarrow} 1$.
\end{proof}

\begin{proof}[\textbf{Proof of Theorem \ref{i1}}]For $j \geq 1$ let $\mathcal{F}_j = \sigma (D_1^{(\vartheta)},\ldots, D_j^{(\vartheta)})$. By the dynamics of the Hoppe tree we have almost surely 
\begin{equation}
\label{EDn}
\E[D_n^{(\vartheta)}|\mathcal{F}_{n-1}]= \frac \vartheta {\vartheta+n-2} (D_1^{(\vartheta)}+1)+ \sum_{i=2}^{n-1} \frac 1 {\vartheta+n-2} (D_i^{(\vartheta)}+1)=1+ \frac 1 {\vartheta+n-2} I_{n-1}^{(\vartheta)}.
\end{equation}
Consequently, $\E[I_n^{(\vartheta)}|\mathcal{F}_{n-1}]=I_{n-1}^{(\vartheta)}+\E[D_n^{(\vartheta)}|\mathcal{F}_{n-1}]=\frac {\vartheta+n-1}{\vartheta+n-2} I_{n-1}^{(\vartheta)} +1$ almost surely.
Therefore, 
\begin{align*}
Z_n^{(\vartheta)}:=\frac 1 {\vartheta+n-1} I_n^{(\vartheta)}-\sum_{i=1}^{n-1}\frac{1}{\vartheta+i}
\end{align*} 
is a zero-mean martingale and $\E[I_n^{(\vartheta)}]=(\vartheta+n-1)\sum_{i=1}^{n-1} \frac 1 {\vartheta+i}$.

The calculations to obtain the expansion for the variance of  $I_n^{(\vartheta)}$ can be done similarly to the calculations in the proof of Theorem \ref{l1}, for details we refer to the master's thesis \cite{Le11}.
\end{proof}

\begin{proof}[\textbf{Proof of Theorem \ref{llipl}}]

To apply a martingale convergence theorem it is sufficient to have a bound on the variance of  the martingale uniformly in $n$. Hence, our expansion of $\V(I_n^{(\vartheta)})$ in Theorem \ref{i1} is sufficient to imply almost sure and $L_2$ convergence of the martingale there, which also applies to the slightly different scaling of  $I_n^{(\vartheta)}$ in Theorem \ref{llipl}.  By our decomposition of the Hoppe tree, see Figure 1, we obtain the recurrence
\begin{align*}
I_n^{(\vartheta)}\stackrel{d}{=} I_{n-N_n}^{(\vartheta)}+ \widetilde{I}^{(1)}_{N_n}+ N_n,
\end{align*}
where $(I_{j}^{(\vartheta)})_{j\ge 1}$, $(\widetilde{I}_{j}^{(1)})_{j\ge 1}$ and $N_n$ are independent and $(\widetilde{I}_{j}^{(1)})_{j\ge 1}$ is distributed as $(I_{j}^{(1)})_{j\ge 1}$. 
For the scaling, 
\begin{align}\label{rec_ipl}
X_n^{(\vartheta)}:=\frac{I_n^{(\vartheta)} - n\log n}{n}
\end{align}
we obtain 
\begin{align}\label{ipl_scaled}
X^{(\vartheta)}_n \stackrel{d}{=} \frac{n-N_n}{n} X_{n-N_n}^{(\vartheta)}+ \frac{N_n}{n}\widetilde{X}^{(1)}_{N_n}+ \frac{1}{n}\left( N_n\log\left(\frac{N_n}{n}\right)+(n-N_n)\log\left(\frac{n-N_n}{n}\right)+N_n\right),
\end{align}
with independence and distributional conditions as in (\ref{rec_ipl}). This suggests that the limit $X^{(\vartheta)}$ of $(X^{(\vartheta)}_n)_{n\ge 1}$ should satisfy the recursive distributional equation
\begin{align}\label{rde}
X^{(\vartheta)}\stackrel{d}{=} (1-B) X^{(\vartheta)}+ B\widetilde{X}^{(1)}+ 
B\log(B)+(1-B)\log(1-B)+B,
\end{align}
where $X^{(\vartheta)}$, $\widetilde{X}^{(1)}$ and $B$ are independent, and $B$ has the beta$(1,\vartheta)$ distribution. Note that $\widetilde{X}^{(1)}$ is the limit distribution of the internal path length of the random recursive tree, that has been obtained by martingale methods by Mahmoud \cite{Ma91} and by the contraction method by Dobrow and Fill \cite{DoFi99}. In particular, in \cite{DoFi99} it is shown that $(X^{(1)}_n)_{n\ge 1}$ converges to its limit $X^{(1)}$ in the minimal $\ell_2$ metric, i.e., weakly and with second moments. This allows us to write the recurrence (\ref{ipl_scaled}) in the form
\begin{align*}
X^{(\vartheta)}_n \stackrel{d}{=} A^{(n)} X_{n-N_n}^{(\vartheta)}+b^{(n)}
\end{align*}
with coefficients
\begin{align*}
A^{(n)}=\frac{n-N_n}{n}, \qquad b^{(n)}=\frac{N_n}{n}\widetilde{X}^{(1)}_{N_n}+ \frac{1}{n}\left( N_n\log\left(\frac{N_n}{n}\right)+(n-N_n)\log\left(\frac{n-N_n}{n}\right)+N_n\right).
\end{align*}
Hence we have convergence of the coefficients to the corresponding quantities in the recursive distributional equation (\ref{rde}) in $\ell_1$, $\ell_2$, in fact in any $\ell_p$, $p\ge 1$. This allows to apply general convergence theorems in the framework of the contraction method, see  R\"osler \cite[Theorem 3]{Ro01} and Neininger and R\"uschendorf \cite[Theorem 4.1]{NeRu04}. In particular, one can first apply Theorem 4.1 in \cite{NeRu04} with the choice of $s=1$ there: This implies convergence in distribution of $X^{(\vartheta)}_n$ to $X^{(\vartheta)}$, where $X^{(\vartheta)}$ is the unique integrable solution of (\ref{rde}), and convergence of the expectations. With this knowledge on the expectation, which, of course, is also covered by our explicit formula for $\E[I_n^{(\vartheta)}]$, one can apply either Theorem 4.1 in \cite{NeRu04} with the choice of $s=2$ or Theorem 3 in \cite{Ro01} to also obtain convergence of the second moments. 

Alternatively to applying the contraction method we could as well use the almost sure convergence of $X^{(\vartheta)}_n$ from the martingale argument together with the almost sure convergence of $N_n/n$ in (\ref{lim_Nn}) to argue that the limit $X^{(\vartheta)}$ satisfies (\ref{rde}).
\end{proof}

\begin{proof}[\textbf{Proof of Theorem \ref{density}}] For the characteristic function $\varphi_\vartheta(t):=\E[\exp(itX^{(\vartheta)})]$ of  $X^{(\vartheta)}$, the recursive distributional equation in Theorem \ref{llipl} implies
\begin{align*}
|\varphi_\vartheta(t)|\le \int_0^1 |\varphi_1(xt)| |\varphi_{\vartheta}((1-x)t)| \vartheta (1-x)^{\vartheta-1}\,dx, \qquad t\in \mathbb{R}.
\end{align*}
We can apply the techniques of Fill and Janson \cite{FJ2000} to show that this relation together with an initial bound on  $|\varphi_\vartheta|$ allows to show that $|\varphi_\vartheta|$ is rapidly decreasing. The details are carried out in the master's thesis \cite{Le11}. Since Fourier transform is an automorphism on the Schwartz space, this implies the assertion.
\end{proof}

\def\cprime{$'$}

\end{document}